\documentclass[12pt]{article}

\usepackage{latexsym}
\usepackage{amsmath}
\usepackage{amssymb}
\usepackage{pstricks}
\usepackage{pst-node}
\usepackage{multido}
\usepackage{calc}
\usepackage{enumerate}
\usepackage{color}
\usepackage{ctable}
\usepackage{tabularx}
\usepackage{xtab}

\makeatletter
\newtheorem{theorem}{Theorem}[section]
\newtheorem{lemma}[theorem]{Lemma}
\newtheorem{conj}[theorem]{Conjecture}

\newtheorem{remark}[theorem]{Remark}

\newtheorem{prop}[theorem]{Proposition}
\newtheorem{defin}[theorem]{Definition}

\def\ps@headings{
 \def\@oddhead{\footnotesize\rm\hfill\runningheadodd\hfill\thepage}
 \def\@evenhead{\footnotesize\rm\thepage\hfill\runningheadeven\hfill}
 \def\@oddfoot{}
 \def\@evenfoot{\@oddfoot}
}

\newcommand{\Prf}{\noindent{\bf Proof}.\quad }

\newcommand{\qed}{\hfill$\Box$}

\makeatother
\def\runningheadeven{odd 2--factored snarks}
\def\runningheadodd{M.Abreu, D.Labbate, R.Rizzi, J.Sheehan}
\title{Odd 2--factored snarks}
\author{{\rm M. Abreu,} \\
\small Dipartimento di Matematica, Informatica ed Economia, Universit\`a della
        Basilicata, \\
      \small C. da Macchia Romana, 85100 Potenza,
         Italy.\\
         \small e-mail: marien.abreu@unibas.it \\
\\
{\rm D. Labbate}, \\
\small Dipartimento di Matematica, Informatica ed Economia, Universit\`a della
        Basilicata, \\
      \small C. da Macchia Romana, 85100 Potenza,
         Italy.\\
        \small e-mail: domenico.labbate@unibas.it \\
\\
{\rm R. Rizzi,} \\
      \small Dipartimento di Informatica,
      Universit\`{a} degli Studi di Verona, \\
      \small Strada le Grazie 15, 37134 Verona, Italy\\
          \small e-mail: romeo.rizzi@univr.it \\
\\
{\rm J. Sheehan,}\\
      \small Department of Mathematical Sciences, King's College,\\
      \small Old Aberdeen AB24 3UE,
         Scotland.\\
        \small e-mail: j.sheehan@maths.abdn.ac.uk }

\date{}

\begin{document}
\maketitle
\pagestyle{headings}

        \begin{abstract}
A {\em snark} is a cubic cyclically $4$--edge connected graph with edge
chromatic number four and girth at least five.
We say that a graph $G$ is {\em odd $2$--factored} if for each $2$--factor F of G each cycle of F is odd.

In this paper, we present a method for constructing odd 2--factored snarks. In particular,
we construct two new odd $2$--factored snarks that disprove a conjecture by some of the authors.
Moreover, we approach the problem of characterizing odd $2$--factored snarks furnishing a partial
characterization of cyclically $4$--edge connected odd $2$--factored snarks. Finally, we pose a new
conjecture regarding odd $2$--factored snarks.
        \end{abstract}

\section{Introduction}

All graphs considered are finite and simple (without loops or
multiple edges). We shall use the term multigraph when multiple
edges are permitted. For definitions and notations not explicitly stated the reader may refer to
\cite{BM}.

A {\em snark} (cf. e.g.~\cite{HS}) is a bridgeless cubic graph with edge chromatic
number four (by Vizing's theorem the edge chromatic number of every cubic
graph is either three or four so a snark corresponds to the special case
of four). In order to avoid trivial cases, snarks are usually assumed to have
girth at least five and not to contain a non--trivial $3$--edge cut (i.e. they are cyclically $4$--edge connected).

Snarks were named after the mysterious and elusive creature in Lewis Caroll's famous poem {\em The Hunting of The Snark} by Martin Gardner in $1976$
\cite{G76}, but it was P. G. Tait in 1880 that initiated the study of snarks, when he proved that the four colour theorem is equivalent to the statement that {\em no snark is planar}~\cite{Ta1880}. The Petersen graph $P$ is the smallest snark and Tutte conjectured that all snarks have Petersen graph minors.
This conjecture was confirmed by Robertson, Seymour and Thomas (cf. \cite{RST}).
Necessarily, snarks are non--hamiltonian.

The importance of the snarks does not only depend on the four colour theorem. Indeed, there are several important open problems
such as the classical cycle double cover conjecture~\cite{Sey79,Sze73}, Fulkerson's conjecture~\cite{Fu71} and Tutte's
5--flow conjecture~\cite{Tu54} for which it is sufficient to prove them for snarks.
Thus, minimal counterexamples to these and other problems must reside, if they exist at all, among the family of snarks.

Snarks play also an important role in characterizing regular graphs with some conditions imposed on their 2--factors.
Recall that a $2$--factor is a $2$--regular spanning subgraph of a graph $G$.

A graph with a $2$--factor is said to be {\em $2$--factor
hamiltonian} if all its $2$--factors are Hamilton cycles, and,
more generally,  {\em $2$--factor isomorphic} if all its
$2$--factors are isomorphic. Examples of such graphs are $K_4$,
$K_5$, $K_{3,3}$, the Heawood graph (which are all $2$--factor
hamiltonian) and the Petersen graph (which is $2$--factor
isomorphic).
Moreover, a {\em pseudo $2$--factor isomorphic graph} is a graphs $G$ with the
property that the parity of the number of cycles in a
$2$--factor is the same for all $2$--factors of $G$. Examples of these graphs are
$K_{3,3}$, the Heawood graph $H_0$ and the Pappus graph $P_0$ (cf. \cite{ADJLS}).
Several papers have addressed the problem of characterizing
families of graphs (particularly regular graphs) which have these
properties directly \cite{Di,FJLS2,AFJLS,AAFJLS,AAFJLS2,FGJ,ADJLS,ALS,ALS2}
or indirectly \cite{FL,L1,L2,FJLS1,L3,AHS,FTV}.
In particular, we have recently pointed out in~\cite{ALS} some relations between snarks and some of these families (cf. Section~\ref{Pre}).

We say that a graph $G$ is {\em odd $2$--factored} (cf.~\cite{ALS}) if for each $2$--factor $F$ of $G$ each cycle of $F$ is odd.
In~\cite{ALS} we have investigated which snarks are odd $2$--factored and we have conjectured that
{\em a snark is odd $2$--factored if and only if $G$ is the
Petersen graph, Blanu\v{s}a~$2$, or a Flower snark $J(t)$, with $t \ge 5$ and odd} (Conjecture~\ref{con1}).

At present, there is no uniform theoretical method for studying snarks and their behaviour. In particular, little is known about the structure of $2$--factors in a given snark.

In this paper,
we present a new method, called {\em bold--gadget dot product}, for constructing odd $2$--factored snarks using the concepts of bold--edges and gadget--pairs
over Isaacs' dot--product~\cite{I75}. This method allows us to construct two new instances of odd $2$--factored snarks of order
$26$ and $34$ that disprove the above conjecture (cf. Conjecture~\ref{con1}). Moreover, we furnish a characterization of bold--edges and gadget--pairs in known odd $2$--factored snarks and we approach the problem of characterizing odd $2$--factored snarks furnishing a partial characterization of cyclically $4$--edge connected odd $2$--factored snarks. Finally, we pose a new conjecture about odd $2$--factored snarks.

\section{Preliminaries}\label{Pre}

Until $1975$ only five snarks were known, then Isaacs~\cite{I75} constructed two infinite families of snarks, one of which is the Flower snark~\cite{I75},
for which in~\cite{ALS} we have used the following definition:


\

\noindent Let $t \geq 5$ be an odd integer. The {\em Flower snark} (cf.~\cite{I75})
$J(t)$ is defined in much the same way as the graph $A(t)$
described in ~\cite{AAFJLS}.

The graph $J(t)$ has vertex set

$$
V(t) = \{h_i, \, u_i, \, v_i, \, w_i \, : \, i=1,2, \ldots , \, t\}
$$

\noindent and edge set

\begin{equation*}
\begin{split}
E(t)& = \{h_iu_i, \, h_iv_i, \, h_iw_i, \, u_iu_{i+1} \,: \, i=1,2, \ldots , \, t\} \\
& \quad \cup\{u_ju_{j+1}, \, v_jv_{j+1}, \, w_jw_{j+1}    \, : \, j=1,2, \ldots , \, t-1\}
  \, \cup \{u_tv_1, \, v_tu_1, \, w_1w_t\}\\
\end{split}
\end{equation*}

For $i=1,2, \ldots , \, t$ we call the subgraph $IC_i$ of $J(t)$
induced by the vertices $\{h_i, \, u_i, \, v_i, \, w_i \}$
the $i^{th}$ {\em interchange} of $J(t)$. The vertices $h_i$
and the edges $\{h_iu_i, \, h_iv_i, \, h_iw_i \}$ are called
respectively the {\em hub} and the {\em spokes} of $IC_i$. The
set of edges $\{ u_iu_{i+1}, \, v_iv_{i+1}, \, w_iw_{i+1}\}$
linking $IC_i$ to $IC_{i+1}$ are said to be the $i^{th}$ {\em link}
$L_i$ of $J(t)$. The edge $u_iu_{i+1} \in L_i$ is called the
$u$--{\em channel of the link}. The subgraph of $J(t)$ induced
by the vertices $\{u_i, \, v_i : \, i=1,2, \ldots , \, t\}$ and
$\{w_i : \, i=1,2, \ldots , \, t\}$ are respectively cycles
of length $2t$ and $t$ and are said to be the {\em base cycles} of $J(t).$

\

The technique used by Isaacs to construct the second infinite family is called a {\em dot product} and it is a consequence of the following:

\begin{lemma}[{\em Parity Lemma}]\cite{I75,Z}\label{paritylemma}
Let $G$ be a cubic graph and let \linebreak $c:E(G) \rightarrow \left\lbrace 1,2,3\right\rbrace $ be a $3$--edge--coloring of $G$.
Then, for every $1$--edge cut $T$ in $G$,
$$ | T \cap c^{-1}(i) | \equiv | T | \; mod \, 2 $$
for each $ i \in \left\lbrace 1,2,3 \right\rbrace. $
\end{lemma}

A \emph{dot product} (see figure below) of two cubic graphs $L$ and $R$, of cyclic--edge--connectivity at least $4$,
denoted by $G= L \cdot R$ is defined as follows~\cite{I75,HS}:

\begin{enumerate}
  \item remove any pair of adjacent vertices $x$ and $y$ from $L$;
  \item remove any two independent edges $ab$ and $cd$ from $R$;
  \item join $\{r,s\}$ to $\{a,b\}$ and $\{t,u\}$ to $\{c,d\}$ or
  $\{r,s\}$ to $\{c,d\}$ and $\{t,u\}$ to $\{a,b\}$, where $N(x)-y=\{r,s\}$
  and $N(y)-x=\{t,u\}$.
\end{enumerate}

\

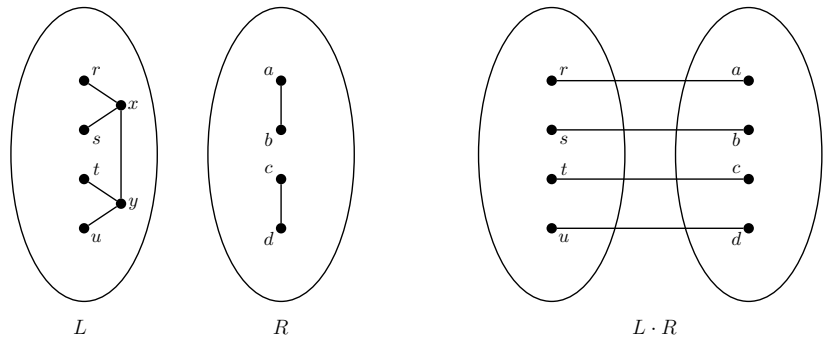
\begin{center}
\resizebox{10cm}{!}{
\begin{minipage}[h]{15cm}
\begin{tabular}{lll}

\begin{pspicture}(-1.5,-2.5)(1.5,2.5)
\psset{unit=1}
\psellipse(0,0)(1.5,3)
\cnode*(0,1.5){3pt}{Nr}\rput(0.25,1.7){\small $r$}
\cnode*(0,0.5){3pt}{Ns}\rput(0.25,0.3){\small $s$}
\cnode*(0,-0.5){3pt}{Nt}\rput(0.25,-0.3){\small $t$}
\cnode*(0,-1.5){3pt}{Nu}\rput(0.25,-1.7){\small $u$}
\cnode*(0.75,1){3pt}{Nx}\rput(1,1){\small $x$}
\cnode*(0.75,-1){3pt}{Ny}\rput(1,-1){\small $y$}

\psellipse(4,0)(1.5,3)
\cnode*(4,1.5){3pt}{Na}\rput(3.75,1.7){\small $a$}
\cnode*(4,0.5){3pt}{Nb}\rput(3.75,0.3){\small $b$}
\cnode*(4,-0.5){3pt}{Nc}\rput(3.75,-0.3){\small $c$}
\cnode*(4,-1.5){3pt}{Nd}\rput(3.75,-1.7){\small $d$}

\ncline{Nx}{Ny}

\ncline{Nr}{Nx}
\ncline{Ns}{Nx}
\ncline{Nt}{Ny}
\ncline{Nu}{Ny}

\ncline{Na}{Nb}
\ncline{Nc}{Nd}

\end{pspicture}

& \hspace{3cm} &

\begin{pspicture}(-1.5,-2.5)(1.5,2.5)
\psset{unit=1}
\psellipse(0,0)(1.5,3)
\cnode*(0,1.5){3pt}{Nr}\rput(0.25,1.7){\small $r$}
\cnode*(0,0.5){3pt}{Ns}\rput(0.25,0.3){\small $s$}
\cnode*(0,-0.5){3pt}{Nt}\rput(0.25,-0.3){\small $t$}
\cnode*(0,-1.5){3pt}{Nu}\rput(0.25,-1.7){\small $u$}

\psellipse(4,0)(1.5,3)
\cnode*(4,1.5){3pt}{Na}\rput(3.75,1.7){\small $a$}
\cnode*(4,0.5){3pt}{Nb}\rput(3.75,0.3){\small $b$}
\cnode*(4,-0.5){3pt}{Nc}\rput(3.75,-0.3){\small $c$}
\cnode*(4,-1.5){3pt}{Nd}\rput(3.75,-1.7){\small $d$}

\ncline{Nr}{Na}
\ncline{Ns}{Nb}
\ncline{Nt}{Nc}
\ncline{Nu}{Nd}

\end{pspicture} \\

\rule{0pt}{1cm}
\hspace{1cm} $L$ \hspace{3.5cm} $R$ & & \hspace{3cm} $L \cdot R$

\end{tabular}
\end{minipage}
} 

\end{center}

\

Note that the dot product allows one to construct graphs of cyclic edge--connectivity exactly $4$.
Moreover, the dot product of the Petersen graph with itself $P \cdot P$ gives rise to two snarks Blanu\v{s}a 1 and Blanu\v{s}a 2.

The Parity Lemma~\ref{paritylemma} allows one to prove the following:

\begin{theorem}\cite{I75,Z}\label{preservesnarks}
Let $L$ and $R$ be snarks. Then the dot product $L \cdot R$ is also a snark.
\end{theorem}

A more general method to construct snarks called {\em superposition} has been introduced by M. Kochol~\cite{K}.
A superposition is performed replacing simultaneously edges and
vertices of a snark by suitable cubic graphs with pendant (or half) edges (called superedges
and supervertices) yielding a new snark.
Superpositions allow one to construct cyclically $k$--edge--connected
snarks with arbitrarily large girth, for $k=4,5,6$.

\

As already mentioned in the Introduction a
graph $G$ is {\em odd $2$--factored}
if for each $2$--factor $F$ of $G$ each cycle of $F$ is odd.

By definition, {\em an odd $2$--factored graph $G$ is pseudo $2$--factor isomorphic}.
Note that,  odd $2$--factoredness is not the same as the {\em oddness} of a (cubic) graph (cf. e.g.\cite{Z}).

\


\begin{lemma}\cite{ALS}\label{snark1}
Let $G$ be a cubic $3$--connected odd $2$--factored  graph then $G$ is a snark.
\end{lemma}

In~\cite{ALS} some of the authors have posed the following:

\

\noindent{\sc Question:} {\em Which snarks are odd $2$--factored?}

\

and we have proved:

\begin{prop}\cite{ALS}\label{snark2}
\begin{enumerate}[(i)]
\item Petersen and Blanu\v{s}a2 are odd 2--factored snarks.

\item The Flower Snark $J(t)$, for odd $t \ge 5$, is odd $2$--factored.
Moreover, $J(t)$ is pseudo $2$--factor isomorphic but
not $2$--factor isomorphic.

\item All other known snarks up to $22$ vertices and
all other named snarks up to $50$ vertices are
not odd $2$--factored.
\end{enumerate}
\end{prop}

Thus it seemed reasonable to pose the following:

\begin{conj}\cite{ALS}\label{con1}
A snark is odd $2$--factored if and only if $G$ is the
Petersen graph, Blanu\v{s}a~$2$, or a Flower snark $J(t)$, with $t \ge 5$ and odd.
\end{conj}

\noindent We disprove Conjecture~\ref{con1} in Section~\ref{counter}.

\

As mentioned above, the Blanu\v{s}a graphs arise as the dot product of the Petersen graph with itself, but one is odd $2$--factored (cf. Proposition~\ref{snark2}(i)) while the other one is not.
In the Petersen graph, which is edge transitive, there are exactly two
kinds of pairs of independent edges. The Blanu\v{s}a snarks are
the result of these two different choices of the pairs of independent
edges in the dot product. We will make use of this property for constructing new odd
$2$--factored snarks in Sections~\ref{constrnewodd} and~\ref{counter}.


\begin{prop}\label{dotnotpreserveodd}
The dot product preserves snarks, but not odd 2–-factored graphs.
\end{prop}

\Prf
It is immediate from Theorem~\ref{preservesnarks} and Proposition~\ref{snark2}(i), (iii).
\qed

\section{A construction of odd $2$--factored snarks}\label{constrnewodd}

We present a general construction of odd $2$--factored snarks performing the dot product on edges with particular properties, called {\em bold--edges} and {\em gadget--pairs} respectively, of two snarks $L$ and $R$.

\

\noindent{\sc Construction: Bold--Gadget Dot Product.}

\

\noindent We construct (new) odd $2$--factored snarks as follows:
\begin{itemize}
\item Take two snarks $L$ and $R$ with bold--edges (cf. Definition~\ref{boldedge}) and gadget--pairs (cf. Definition~\ref{gadgetpair}), respectively;
\item Choose a bold--edge $xy$ in $L$;
\item Choose a gadget--pair $f$, $g$ in $R$;
\item Perform a dot product $L \cdot R$ using these edges;
\item Obtain a new odd 2--factored snark (cf. Theorem~\ref{OddDot}).
\end{itemize}

Note that in what follows the existence of a 2--factor in a snark is guaranteed since they are bridgeless by definition.

\begin{defin}\label{boldedge}
Let $L$ be a snark. A {\em bold--edge} is an edge $e=xy \in L$ such that the following conditions hold:

\begin{enumerate}[(i)]
  \item All $2$--factors of $L-x$ and of $L-y$ are odd;
  \item all $2$--factors of $L$ containing $xy$ are odd;
  \item all $2$--factors of $L$ avoiding $xy$ are odd.
\end{enumerate}
\end{defin}

Note that not all snarks contain bold--edges (cf. Proposition~\ref{P34}, Lemma~\ref{FlowerNoBold}). Furthermore, conditions $(ii)$ and $(iii)$ are trivially satisfied if $L$ is odd $2$--factored.

\begin{lemma}\label{PetBold}
The edges of the Petersen graph $P_{10}$ are all bold--edges.
\end{lemma}

\Prf
Since $P_{10}$ is hypohamiltonian (i.e. $P_{10}-v$ is hamiltonian, for each $v \in V(P_{10})$) and moreover, for every $v \in P_{10}$, all $2$--factors of $P_{10}-v$ are hamiltonian , condition $(i)$ holds. The other two conditions are satisfied since $P_{10}$ is odd-$2$--factored.
\qed

\begin{defin}\label{gadgetpair}
Let $R$ be a snark. A pair of independent edges $f=ab$ and $g=cd$ is called a {\em gadget--pair} if the following conditions hold:

\begin{enumerate}[(i)]
  \item There are no $2$--factors of $R$ avoiding both $f,g$;
  \item all $2$--factors of $R$ containing exactly one element of $\{f,g\}$ are odd;

  \item all $2$--factors of $R$ containing both $f$ and $g$ are odd. Moreover, $f$ and $g$ belong to different cycles in each such factor.

  \item all $2$--factors of $(R-\{f,g\}) \cup \{ac,ad,bc,bd\}$ containing exactly one
  element of $\{ac,ad,bc,bd\}$, are such that the cycle containing the new edge is even and all other cycles are odd.
\end{enumerate}
\end{defin}

Note that, finding gadget--pairs in a snark
is not an easy task and, in general, not all snarks contain gadget--pairs (cf. Lemma~\ref{gadgetsearch}).

\begin{figure}[h]
\begin{center}
\resizebox{3cm}{!}{
\begin{minipage}[h]{5cm}
\begin{pspicture}(-1.5,-1.5)(1.5,1.5)
\psset{unit=0.75}
\psset{linewidth=0.5pt}
\SpecialCoor
\degrees[5]
\cnode*(2.5;1.25){3pt}{P11}
\cnode*(2.5;2.25){3pt}{P12}
\cnode*(2.5;3.25){3pt}{P13}
\cnode*(2.5;4.25){3pt}{P14}
\cnode*(2.5;5.25){3pt}{P15}
\cnode*(1.25;1.25){3pt}{P21}
\cnode*(1.25;2.25){3pt}{P22}
\cnode*(1.25;3.25){3pt}{P23}
\cnode*(1.25;4.25){3pt}{P24}
\cnode*(1.25;5.25){3pt}{P25}
\ncline{P11}{P12}
\ncline{P12}{P13}
\ncline{P14}{P15}
\ncline{P15}{P11}
\ncline{P12}{P22}
\ncline{P13}{P23}
\ncline{P14}{P24}
\ncline{P15}{P25}
\ncline{P21}{P24}
\ncline{P24}{P22}
\ncline{P25}{P23}
\ncline{P23}{P21}
\psset{linewidth=1pt,linestyle=dashed}
\ncline{P11}{P21}
\ncline{P22}{P25}
\ncline{P13}{P14}
\end{pspicture}
\end{minipage}
}
\caption{Any pair of the dashed edges is a gadget--pair in $P_{10}$}
\label{PetSpc}
\end{center}
\end{figure}
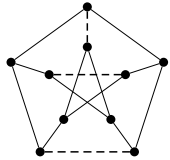

Let $H:=\{x_1y_1,x_2y_2,x_3y_3\}$ be the two horizontal
edges and the vertical edge respectively (in the pentagon--pentagram representation) of $P_{10}$ (cf. Figure~\ref{PetSpc}).

It is easy to prove the following properties:

\begin{lemma}\label{petersenproperties}
Let $P_{10}$ be the Petersen graph and $H:=\{x_1y_1,x_2y_2,x_3y_3\}$ be as above.

\begin{enumerate}[(i)]
\item The graph $P_{10}-H$ is bipartite.
\item The graph $P_{10}-\{f,g\}$ has no 2--factors, for any distinct $f,g \in  H$.
\end{enumerate}

\end{lemma}

\begin{lemma}\label{PetGadg}
Any pair of distinct edges $f,g$ in the set $H$ of $P_{10}$ is a gadget--pair.
\end{lemma}

\Prf
It can be easily checked that the edges of $H$ in $P_{10}$ have the property that any $2$--factor of $P_{10}$ contains exactly
two of them. Moreover, they belong to different cycles of the $2$--factor. Indeed, for any two edges $f$ and $g$ in $H$ their endvertices are all at distance $2$ in $P_{10}$.  Thus the shortest cycle containing both $f$ and $g$ has length $6$. Since all $2$--factors of $P_{10}$ contain two $5$--cycles, in any $2$--factor of $P_{10}$ containing both $f$ and $g$, these edges are contained in different cycles. Hence, conditions $(i)$--$(iii)$ follow from the above reasoning and the odd $2$--factoredness of $P_{10}$.

Condition $(iv)$ can also be easily checked and moreover,
any $2$--factor of $P_{10}-\{f,g\}+\{x_ix_j\}$ (or  $P_{10}-\{f,g\}+\{x_iy_j\}$ or $P_{10}-\{f,g\}+\{y_iy_j\}$), for $i \ne j$,
containing the new edge is hamiltonian, hence even (and obviously there are no other cycles in these 2--factors).
\qed

\

In the next lemma we recall a well known property of edge--cuts:

\begin{lemma}\label{evennumberofedgesinT}
Let $G$ be a connected graph and let $S$ be a set of edges such that $G-S$ is disconnected, but $G-S'$ is not disconnected,
for any proper subset $S'$ of $S$. Then, for any cycle $C$ of $G$, $E(C) \cap S$ is even.
\end{lemma}



Recall that the {\em length} of a cycle $C$ is denoted by $|C|$. The following theorem allows us to construct new odd $2$--factored snarks.

\begin{theorem} \label{OddDot}
Let $xy$ be a bold--edge in a snark $L$ and let $\{ab,cd\}$ be a gadget--pair in a snark $R$.
Then $L \cdot R$ is an odd $2$--factored snark.
\end{theorem}

\Prf
Denote $e:=xy$, $f:=ab$, $g:=cd$, $N(x)-y:=\{r,s\}$, $N(y)-x:=\{t,u\}$ and $T:=\{ra,sb,tc,ud\}$ the 4--edge cut obtained performing the dot product $L \cdot R$.

Let $F$ be a 2--factor of $L \cdot R$ then $F$ contains an even number of edges of $T$ by Lemma~\ref{evennumberofedgesinT}.

We distinguish three cases according to the number of edges of $T$ in $F$:

\

\noindent {\sc Case 1.} $F$ contains no edges of $T$.

In this case it is immediate to check that a subset of the cycles of $F$ forms a $2$--factor of $R-\{f,g\}$, contradicting Definition~\ref{gadgetpair}(i).
Thus there are no 2--factors of $L \cdot R$ avoiding $T$.

\

\noindent {\sc Case 2.} $F$ contains exactly two edges $e_1$ and $e_2$ of $T$.

We want to prove that all cycles of $F$ are odd. We distinguish two subcases.

\noindent {\sc Case 2.1.} The endvertices of $e_1$ and $e_2$ in $R$ are both endvertices of either $f$ or $g$.

W.l.g. we may assume that $e_1=ra$ and $e_2=sb$. Let $F_1:=F \cap L$ and let $F_1':=F_1 \cup \{rx,sx\}$. Then $F_1'$ is a 2--factor of $L-y$.
Analogously, let $F_2:=F \cap R$ and let $F_2':=F_2 \cup \{f\}$. Then $F_2'$ is a 2--factor of $R$ containing $f$ and avoiding $g$. Let $C_x$ be the cycle of $F_1'$ containing $x$. Then $|C_x|$ is odd by Definition~\ref{boldedge}(i). Similarly, let $C_f$ be the cycle of $F_2'$ containing $f$. Then $|C_f|$ is odd by Definition~\ref{gadgetpair}(ii). Thus, the cycle $C$ of $F$ containing $e_1$ and $e_2$ has $|C|=|C_x|-2+|C_f|-1+2=|C_x|+|C_f|-1$ which is odd. Finally,all other cycles of $F$ are odd by Definition~\ref{boldedge}(i) and Definition~\ref{gadgetpair}(ii).

\noindent {\sc Case 2.2.} The endvertices of $e_1$ and $e_2$ in $R$ lie one in $f$ and the other in $g$.

W.l.g. we may assume that $e_1=ra$ and $e_2=tc$. Let $F_1:=F \cap L$ and let $F_1':=F_1 \cup \{rx,xy,yt\}$. Then $F_1'$ is a 2--factor of $L$ containing $xy$. Analogously, let $F_2:=F \cap R$ and let $F_2':=F_2 \cup \{ac\}$. Let $S:=\{ac,ad,bc,bd\}$ be a set of new edges and consider the graph $R':=R-\{f,g\} \cup S$. Then $F_2'$ is a 2--factor of $R'$ containing only $ac$ of $S$ by construction. Let $C_{xy}$ be the cycle of $F_1'$ containing $xy$. Then $|C_{xy}|$ is odd by Definition~\ref{boldedge}(ii). Similarly, let $C_{ac}$ be the cycle of $F_2'$ containing $ac$. Then $|C_{ac}|$ is even by Definition~\ref{gadgetpair}(iv). Thus, the cycle $C$ of $F$ containing $e_1$ and $e_2$ has $|C|=|C_{xy}|-3+|C_{ac}|-1+2=|C_{xy}|+|C_{ac}|-2$ which is odd. Finally, all other cycles of $F$ are odd by Definition~\ref{boldedge}(ii) and Definition~\ref{gadgetpair}(iv).

\

\noindent {\sc Case 3.} $F$ contains all the four edges of $T$.

Again we want to prove that all cycles of $F$ have odd length. Let $F_1:=F \cap L$, $F_2:=F \cap R$, $F_1':=F_1\cup \{rx,sx,ty,uy\}$ and $F_2':=F_2 \cup \{ab,cd\}$. Note that $F_1'$ is a $2$--factor of $L$ avoiding $xy$ and that $F_2'$ is a $2$--factor of $R$ containing both $f$ and $g$.
Let $C_x$ and $C_y$ be the cycles of $F_1'$ containing $x$ and $y$, respectively. If $C_x=C_y$ then we denote such a cycle by $C_{xy}$.
Analogously, let $C_f$ and $C_g$ be the cycles of $F_2'$ containing $f$ and $g$, respectively. Note that $C_f$ and $C_g$ are always distinct by Definition~\ref{gadgetpair}(iii).

In order to compute the parity of the length of the cycles of $F$ containing $T$,
we need to analyze all possible combinations of paths in $F$ between the vertices $\{r,s,t,u\}$ and between the vertices $\{a,b,c,d\}$ of $L \cdot R$.
It is easy to check that we have five different cases (the others being equivalent to some of these five) but three of them are ruled out by Definition~\ref{gadgetpair}(iii), since they have $C_f = C_g$  (cf. Figure~\ref{Subcases3}).

\begin{figure}[h]

\begin{center}
\begin{tabular}{cc}

\resizebox{1.5cm}{!}{
\begin{minipage}[h]{7.5cm}
\begin{pspicture}(-3.75,-6.25)(3.75,6.25)
\psset{unit=1}
\psellipse(0,0)(3.75,7.5)
\cnode*(0,3.75){7.5pt}{Nr}\rput(0.625,4.25){\Huge $r$}
\cnode*(0,1.25){7.5pt}{Ns}\rput(0.625,0.75){\Huge $s$}
\cnode*(0,-1.25){7.5pt}{Nt}\rput(0.625,-0.75){\Huge $t$}
\cnode*(0,-3.75){7.5pt}{Nu}\rput(0.625,-4.25){\Huge $u$}

\psellipse(10,0)(3.75,7.5)
\cnode*(10,3.75){7.5pt}{Na}\rput(9.375,4.25){\Huge $a$}
\cnode*(10,1.25){7.5pt}{Nb}\rput(9.375,0.75){\Huge $b$}
\cnode*(10,-1.25){7.5pt}{Nc}\rput(9.375,-0.75){\Huge $c$}
\cnode*(10,-3.75){7.5pt}{Nd}\rput(9.375,-4.25){\Huge $d$}

\ncline{Nr}{Na}
\ncline{Ns}{Nb}
\ncline{Nt}{Nc}
\ncline{Nu}{Nd}

\pscurve(0,3.75)(-1.5,2.5)(-0.75,1.25)(-1.5,0)(-0.75,-1.25)(-1.5,-2.5)(0,-3.75)
\pscurve(0,1.25)(-0.5,0.9375)(-0.25,0.625)(-0.5,0)(-0.25,-0.625)(-0.5,-0.9375)(0,-1.25)

\pscurve(10,1.25)(10.75,1.5625)(10.25,1.875)(10.75,2.1875)(10.25,2.5)(10.75,2.8125)(10.25,3.125)(10.75,3.4375)(10,3.75)
\pscurve(10,-1.25)(10.75,-1.5625)(10.25,-1.875)(10.75,-2.1875)(10.25,-2.5)(10.75,-2.8125)(10.25,-3.125)(10.75,-3.4375)(10,-3.75)

\end{pspicture}
\end{minipage}
}

&

\hspace{0.25cm}
\resizebox{1.5cm}{!}{
\begin{minipage}[h]{7.5cm}
\begin{pspicture}(-3.75,-6.25)(3.75,6.25)
\psset{unit=1}
\psellipse(0,0)(3.75,7.5)
\cnode*(0,3.75){7.5pt}{Nr}\rput(0.625,4.25){\Huge $r$}
\cnode*(0,1.25){7.5pt}{Ns}\rput(0.625,0.75){\Huge $s$}
\cnode*(0,-1.25){7.5pt}{Nt}\rput(0.625,-0.75){\Huge $t$}
\cnode*(0,-3.75){7.5pt}{Nu}\rput(0.625,-4.25){\Huge $u$}

\psellipse(10,0)(3.75,7.5)
\cnode*(10,3.75){7.5pt}{Na}\rput(9.375,4.25){\Huge $a$}
\cnode*(10,1.25){7.5pt}{Nb}\rput(9.375,0.75){\Huge $b$}
\cnode*(10,-1.25){7.5pt}{Nc}\rput(9.375,-0.75){\Huge $c$}
\cnode*(10,-3.75){7.5pt}{Nd}\rput(9.375,-4.25){\Huge $d$}

\ncline{Nr}{Na}
\ncline{Ns}{Nb}
\ncline{Nt}{Nc}
\ncline{Nu}{Nd}

\pscurve(0,1.25)(-0.75,1.5625)(-0.25,1.875)(-0.75,2.1875)(-0.25,2.5)(-0.75,2.8125)(-0.25,3.125)(-0.75,3.4375)(0,3.75)
\pscurve(0,-1.25)(-0.75,-1.5625)(-0.25,-1.875)(-0.75,-2.1875)(-0.25,-2.5)(-0.75,-2.8125)(-0.25,-3.125)(-0.75,-3.4375)(0,-3.75)

\pscurve(10,1.25)(10.75,1.5625)(10.25,1.875)(10.75,2.1875)(10.25,2.5)(10.75,2.8125)(10.25,3.125)(10.75,3.4375)(10,3.75)
\pscurve(10,-1.25)(10.75,-1.5625)(10.25,-1.875)(10.75,-2.1875)(10.25,-2.5)(10.75,-2.8125)(10.25,-3.125)(10.75,-3.4375)(10,-3.75)

\end{pspicture}
\end{minipage}
} \\

\rule{0pt}{0.75cm}
\hspace{2cm}{\sc Case 3.1} &

\hspace{2cm}
{\sc Case 3.2}

\end{tabular}
\end{center}

\begin{center}
\begin{tabular}{ccc}

\rule{0pt}{2cm}
\resizebox{1.5cm}{!}{
\begin{minipage}[h]{7.5cm}
\begin{pspicture}(-3.75,-6.25)(3.75,6.25)
\psset{unit=1}
\psellipse(0,0)(3.75,7.5)
\cnode*(0,3.75){7.5pt}{Nr}\rput(0.625,4.25){\Huge $r$}
\cnode*(0,1.25){7.5pt}{Ns}\rput(0.625,0.75){\Huge $s$}
\cnode*(0,-1.25){7.5pt}{Nt}\rput(0.625,-0.75){\Huge $t$}
\cnode*(0,-3.75){7.5pt}{Nu}\rput(0.625,-4.25){\Huge $u$}

\psellipse(10,0)(3.75,7.5)
\cnode*(10,3.75){7.5pt}{Na}\rput(9.375,4.25){\Huge $a$}
\cnode*(10,1.25){7.5pt}{Nb}\rput(9.375,0.75){\Huge $b$}
\cnode*(10,-1.25){7.5pt}{Nc}\rput(9.375,-0.75){\Huge $c$}
\cnode*(10,-3.75){7.5pt}{Nd}\rput(9.375,-4.25){\Huge $d$}

\ncline{Nr}{Na}
\ncline{Ns}{Nb}
\ncline{Nt}{Nc}
\ncline{Nu}{Nd}

\pscurve(0,1.25)(-0.75,1.5625)(-0.25,1.875)(-0.75,2.1875)(-0.25,2.5)(-0.75,2.8125)(-0.25,3.125)(-0.75,3.4375)(0,3.75)
\pscurve(0,-1.25)(-0.75,-1.5625)(-0.25,-1.875)(-0.75,-2.1875)(-0.25,-2.5)(-0.75,-2.8125)(-0.25,-3.125)(-0.75,-3.4375)(0,-3.75)

\pscurve(10,3.75)(11.5,2.5)(10.75,1.25)(11.5,0)(10.75,-1.25)(11.5,-2.5)(10,-3.75)
\pscurve(10,1.25)(10.5,0.9375)(10.25,0.625)(10.5,0)(10.25,-0.625)(10.5,-0.9375)(10,-1.25)

\end{pspicture}
\end{minipage}
}

&

\hspace{2cm}
\resizebox{1.5cm}{!}{
\begin{minipage}[h]{7.5cm}
\begin{pspicture}(-3.75,-6.25)(3.75,6.25)
\psset{unit=1}
\psellipse(0,0)(3.75,7.5)
\cnode*(0,3.75){7.5pt}{Nr}\rput(0.625,4.25){\Huge $r$}
\cnode*(0,1.25){7.5pt}{Ns}\rput(0.625,0.75){\Huge $s$}
\cnode*(0,-1.25){7.5pt}{Nt}\rput(0.625,-0.75){\Huge $t$}
\cnode*(0,-3.75){7.5pt}{Nu}\rput(0.625,-4.25){\Huge $u$}

\psellipse(10,0)(3.75,7.5)
\cnode*(10,3.75){7.5pt}{Na}\rput(9.375,4.25){\Huge $a$}
\cnode*(10,1.25){7.5pt}{Nb}\rput(9.375,0.75){\Huge $b$}
\cnode*(10,-1.25){7.5pt}{Nc}\rput(9.375,-0.75){\Huge $c$}
\cnode*(10,-3.75){7.5pt}{Nd}\rput(9.375,-4.25){\Huge $d$}

\ncline{Nr}{Na}
\ncline{Ns}{Nb}
\ncline{Nt}{Nc}
\ncline{Nu}{Nd}

\pscurve(0,3.75)(-1.5,2.5)(-0.75,1.25)(-1.5,0)(-0.75,-1.25)(-1.5,-2.5)(0,-3.75)
\pscurve(0,1.25)(-0.5,0.9375)(-0.25,0.625)(-0.5,0)(-0.25,-0.625)(-0.5,-0.9375)(0,-1.25)

\pscurve(10,3.75)(11.5,2.5)(10.75,1.25)(11.5,0)(10.75,-1.25)(11.5,-2.5)(10,-3.75)
\pscurve(10,1.25)(10.5,0.9375)(10.25,0.625)(10.5,0)(10.25,-0.625)(10.5,-0.9375)(10,-1.25)

\end{pspicture}
\end{minipage}
}

&

\hspace{1cm}
\resizebox{1.5cm}{!}{
\begin{minipage}[h]{7.5cm}
\begin{pspicture}(-3.75,-6.25)(3.75,6.25)
\psset{unit=1}
\psellipse(0,0)(3.75,7.5)
\cnode*(0,3.75){7.5pt}{Nr}\rput(0.625,4.25){\Huge $r$}
\cnode*(0,1.25){7.5pt}{Ns}\rput(0.625,0.75){\Huge $s$}
\cnode*(0,-1.25){7.5pt}{Nt}\rput(0.625,-0.75){\Huge $t$}
\cnode*(0,-3.75){7.5pt}{Nu}\rput(0.625,-4.25){\Huge $u$}

\psellipse(10,0)(3.75,7.5)
\cnode*(10,3.75){7.5pt}{Na}\rput(9.375,4.25){\Huge $a$}
\cnode*(10,1.25){7.5pt}{Nb}\rput(9.375,0.75){\Huge $b$}
\cnode*(10,-1.25){7.5pt}{Nc}\rput(9.375,-0.75){\Huge $c$}
\cnode*(10,-3.75){7.5pt}{Nd}\rput(9.375,-4.25){\Huge $d$}

\ncline{Nr}{Na}
\ncline{Ns}{Nb}
\ncline{Nt}{Nc}
\ncline{Nu}{Nd}

\pscurve(0,1.25)(-1.5,0.625)(-0.75,0)(-1.5,-0.625)(-0.75,-1.25)(-1.5,-1.875)(-0.75,-2.5)(-1.5,-3.125)(0,-3.75)
\pscurve(0,-1.25)(-0.75,-0.5)(-0.25,0)(-0.75,1.5)(-0.25,2.5)(-0.75,3.25)(0,3.75)

\pscurve(10,3.75)(11.5,2.5)(10.75,1.25)(11.5,0)(10.75,-1.25)(11.5,-2.5)(10,-3.75)
\pscurve(10,1.25)(10.5,0.9375)(10.25,0.625)(10.5,0)(10.25,-0.625)(10.5,-0.9375)(10,-1.25)

\end{pspicture}
\end{minipage}
} \\

\multicolumn{3}{c}{\rule{0pt}{1cm} These three cases can be ruled out since they give rise to} \\

\multicolumn{3}{c}{$C_f=C_g$ in $F_2'$, impossible by Definition~\ref{gadgetpair}(iii)}

\end{tabular}
\end{center}

  \caption{Subcases of \sc Case 3 }\label{Subcases3}
\end{figure}

The two remaining subcases are:

\

\noindent {\sc Case 3.1.} All edges of $T$ lie in a cycle $C$ of $F$ such that $C_x=C_y=C_{xy}$ and $C_f \neq C_g$.

In this case $|C|=|C_{xy}|-4+|C_f|-1+|C_g|-1+4$$=|C_{xy}|+|C_f|+|C_g|-2$ which is odd by Definition~\ref{boldedge}(iii) and Definition~\ref{gadgetpair}(iii).

\noindent {\sc Case 3.2.} The edges of $T$ are contained in two distinct cycles $C_1$ and $C_2$ of $F$ such that  $C_x \neq C_y$ and  $C_f \neq C_g$.

In this case $|C_1|=|C_x|-2+|C_f|-1+2$$=|C_x|+|C_f|-1$ and $|C_2|=|C_y|-2+|C_g|-1+2$$=|C_y|+|C_g|-1$ which are odd by Definition~\ref{boldedge}(iii) and Definition~\ref{gadgetpair}(iii).

All the remaining cycles of $F$ are odd by Definition~\ref{boldedge}(iii) and Definition~\ref{gadgetpair}(iii).

Thus, the resulting graph is odd $2$--factored, hence a snark by Lemma~\ref{snark1}.
\qed

\

\noindent{\sc Construction of $P_{18}$}

\

Recall that, the Blanu\v{s}a2 snark is odd $2$--factored (cf.~\cite{ALS} and Proposition~\ref{snark2}).
We can obtain the same result taking two copies $L,R$ of the Petersen graph $P_{10}$, in the first one choosing any edge as a bold--edge (by Lemma~\ref{PetBold}) and in the second a gadget--pair as in Lemma~\ref{PetGadg}. The resulting graph, obtained
as the dot product $L \cdot R$, denoted by $P_{18}$, is odd $2$--factored by Theorem~\ref{OddDot} (cf. Figure~\ref{P18Fig}) and it is isomorphic to the Blanu\v{s}a2 snark.

\

Let $H:=\{e_1,e_2,e_3\}$ be the two horizontal
edges and the vertical edge respectively (in the pentagon--pentagram representation) of the Petersen
graph $P_{10}$, as in Figure~\ref{PetSpc}. Let $L$ and $R$ be two copied of $P_{10}$.
Choose $e_1=xy$ as the bold--edge in $L$ and $f,g \in H$ as the gadget--pair in $R$.
Moreover, let  $L_0:=L-\{x,y\}$ and $R_0:=R-\{f,g\}$ be the {\em $4$--poles} represented as follows:

\begin{figure}[h]
\begin{center}
\resizebox{6cm}{!}{
\begin{minipage}[h]{12cm}
\begin{pspicture}(-1,-3)(8,4)
\psset{unit=1}
\cnode*(-2,0){3pt}{N1}
\cnode*(-1,-1){3pt}{N2}
\cnode*(0,-2){3pt}{N3}
\cnode*(-1,1){3pt}{N4}
\cnode*(0,2){3pt}{N5}
\cnode*(1,1){3pt}{N6}
\cnode*(2,0){3pt}{N7}
\cnode*(1,-1){3pt}{N8}
\cnode*(12,0){3pt}{N9}
\cnode*(11,-1){3pt}{N10}
\cnode*(10,-2){3pt}{N11}
\cnode*(10.5,-0.5){3pt}{N12}
\cnode*(10.5,0.5){3pt}{N13}
\cnode*(11,1){3pt}{N14}
\cnode*(10,2){3pt}{N15}
\cnode*(9,1){3pt}{N16}
\cnode*(8,0){3pt}{N17}
\cnode*(9,-1){3pt}{N18}
\ncline{N1}{N2}
\ncline{N1}{N4}
\ncline{N2}{N3}
\ncline{N2}{N6}\rput*(0.65,0.05){\small $e_2$}
\ncline{N3}{N8}
\ncline{N4}{N5}
\ncline{N4}{N8}\rput*(-0.65,0.05){\small $e_3$}
\ncline{N5}{N6}
\ncline{N6}{N7}
\ncline{N7}{N8}
\ncline{N9}{N10}
\ncline{N9}{N14}
\ncline{N10}{N11}
\ncline{N10}{N12}
\ncline{N11}{N18}
\ncline{N12}{N13}
\ncline{N12}{N16}
\ncline{N13}{N14}
\ncline{N13}{N18}
\ncline{N14}{N15}
\ncline{N15}{N16}
\ncline{N16}{N17}
\ncline{N17}{N18}
\cnode*(3,3){0.5pt}{N1p}
\ncangle[angleA=90,angleB=180,arm=1cm,linearc=.15]{N1}{N1p}
\psline(0,-2)(3,-2)
\psline(0,2)(3,2)
\psline(2,0)(3,0)
\cnode*(7,3){0.5pt}{N9p}
\ncangle[angleA=90,angleB=0,arm=1cm,linearc=.15]{N9}{N9p}
\psline(10,-2)(7,-2)
\psline(10,2)(7,2)
\psline(8,0)(7,0)
\rput(0,-3){\Huge $L_0 := P_{10}-\{x,y\}$}
\rput(10,-3){\Huge $R_0 := P_{10}-\{f,g\}$}
\end{pspicture}
\end{minipage}
}
\end{center}
\caption{4--poles from Petersen}\label{4-polesfig}
\end{figure}
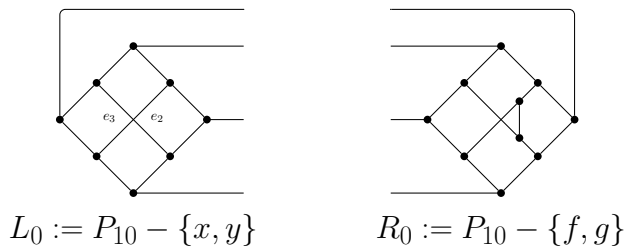

\noindent Performing the dot--product $L \cdot R$ we obtain the Blanu\v{s}a2 snark $P_{18}$ (Figure~\ref{P18Fig}).

\

\begin{figure}[h]
\begin{center}
\resizebox{6cm}{!}{
\begin{minipage}[h]{12cm}
\begin{pspicture}(-3,-2)(4,2)
\psset{unit=1}
\cnode*(-2,0){3pt}{N1}\rput(-2.25,0){\small $1$}
\cnode*(-1,-1){3pt}{N2}\rput(-1.25,-1.25){\small $2$}
\cnode*(0,-2){3pt}{N3}\rput(0,-2.25){\small $3$}
\cnode*(-1,1){3pt}{N4}\rput(-1.25,1.25){\small $4$}
\cnode*(0,2){3pt}{N5}\rput(0,2.25){\small $5$}
\cnode*(1,1){3pt}{N6}\rput(1.15,1.25){\small $6$}
\cnode*(2,0){3pt}{N7}\rput(2,-0.35){\small $7$}
\cnode*(1,-1){3pt}{N8}\rput(1.15,-1.25){\small $8$}
\cnode*(8,0){3pt}{N9}\rput(8,-0.35){\small $9$}
\cnode*(7,-1){3pt}{N10}\rput(7.15,-1.25){\small $10$}
\cnode*(6,-2){3pt}{N11}\rput(6,-2.25){\small $11$}
\cnode*(6.5,-0.5){3pt}{N12}\rput(6.15,-0.75){\small $12$}
\cnode*(6.5,0.5){3pt}{N13}\rput(6.15,0.75){\small $13$}
\cnode*(7,1){3pt}{N14}\rput(7.15,1.25){\small $14$}
\cnode*(6,2){3pt}{N15}\rput(6,2.25){\small $15$}
\cnode*(5,1){3pt}{N16}\rput(4.75,1.25){\small $16$}
\cnode*(4,0){3pt}{N17}\rput(4,-0.35){\small $17$}
\cnode*(5,-1){3pt}{N18}\rput(4.75,-1.25){\small $18$}
\ncline{N1}{N2}
\ncline{N1}{N4}
\ncline{N2}{N3}
\ncline[linewidth=2pt]{N2}{N6}\rput*(0.65,0.05){\small $e_2$}
\ncline{N3}{N8}
\ncline{N4}{N5}
\ncline[linewidth=2pt]{N4}{N8}\rput*(-0.65,0.05){\small $e_3$}
\ncline{N5}{N6}
\ncline{N6}{N7}
\ncline{N7}{N8}
\ncline{N9}{N10}
\ncline{N9}{N14}
\ncline{N10}{N11}
\ncline{N10}{N12}
\ncline{N11}{N18}
\ncline{N12}{N13}
\ncline{N12}{N16}
\ncline{N13}{N14}
\ncline{N13}{N18}
\ncline{N14}{N15}
\ncline{N15}{N16}
\ncline{N16}{N17}
\ncline{N17}{N18}
\ncangles[angleA=90,angleB=90,arm=3cm,linearc=.15]{N1}{N9}
\ncline{N5}{N15}
\ncline{N7}{N17}
\ncline{N3}{N11}
\end{pspicture}
\end{minipage}
}
\end{center}

\caption{$P_{18}$}
\label{P18Fig}
\end{figure}


\begin{lemma}\label{P18}
Under the above hypothesis, the only bold--edges of $P_{18}$ are
those edges, say $e_2$ and $e_3$, identified with the edges $e_2$ and $e_3$ of $L$ (cf. Figure~\ref{P18Fig}).
\end{lemma}

\Prf
Fix the labelling on $V(P_{18})$ as in Figure~\ref{P18Fig}. To find all possible bold--edges in $P_{18}$,
we only need to verify Definition~\ref{boldedge}$(i)$, since $P_{18}$ is odd $2$--factored.

To this purpose, we have implemented a program, with the software package MAGMA~\cite{BCP}, and computed that
the graph $P_{18}$ has the dihedral group $D_4$ as automorphism group, its edge--orbits are six and its vertex--orbits are five.
For each representative $v$ of the vertex--orbits, we have determined all the $2$--factors of  $P_{18}-v$
(computing the determinant of the variable adjacency matrix of $G$~\cite{Ha62}).
The only vertex, for which $P_{18}-v$ has only odd $2$--factors, is $v=2$ (lying in the vertex--orbit $\{2,4,6,8\}$).
Hence, the only bold--edges in $P_{18}$ are $e_2,e_3$, since there is an edge--orbit $E_0:=\{(2,6), (4,8)\}$ of $P_{18}$, and its edges correspond to $e_2,e_3$ (c.f. Figure~\ref{P18Fig}).
\qed

\section{Counterexamples to Conjecture~\ref{con1}: Constructions}\label{counter}

We construct two new examples of odd $2$--factored snarks of order $26$ and $34$, denoted respectively as $P_{26}$ and $P_{34}$, and starting from the Petersen and the Blanu\v{s}a2 snarks applying iteratively the method described in Section~\ref{constrnewodd}. These two examples disprove Conjecture~\ref{con1}. Moreover, we investigate the structure of the snarks obtained with this method
computing their bold--edges and gadget pairs.

\


\noindent {\sc Construction of $P_{26}$}

\

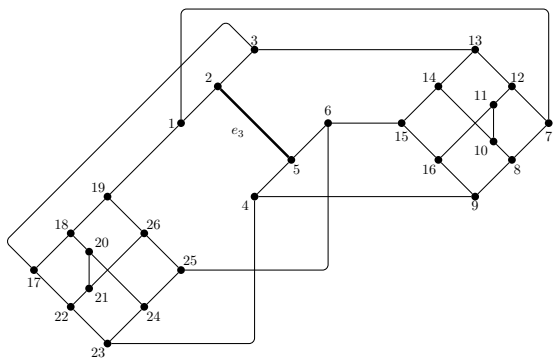
\begin{figure}[h]
\begin{center}
\resizebox{6cm}{!}{
\begin{minipage}[h]{12cm}
\begin{pspicture}(-4,-6)(6,4)
\psset{unit=1}
\cnode*(-2,0){3pt}{N1}\rput(-2.25,0){\small $1$}
\cnode*(0,-2){3pt}{N4}\rput(-0.25,-2.25){\small $4$}
\cnode*(-1,1){3pt}{N2}\rput(-1.25,1.25){\small $2$}
\cnode*(0,2){3pt}{N3}\rput(0,2.25){\small $3$}
\cnode*(2,0){3pt}{N6}\rput(2,0.35){\small $6$}
\cnode*(1,-1){3pt}{N5}\rput(1.15,-1.25){\small $5$}
\cnode*(8,0){3pt}{N7}\rput(8,-0.35){\small $7$}
\cnode*(7,-1){3pt}{N8}\rput(7.15,-1.25){\small $8$}
\cnode*(6,-2){3pt}{N9}\rput(6,-2.25){\small $9$}
\cnode*(6.5,-0.5){3pt}{N10}\rput(6.15,-0.75){\small $10$}
\cnode*(6.5,0.5){3pt}{N11}\rput(6.15,0.75){\small $11$}
\cnode*(7,1){3pt}{N12}\rput(7.15,1.25){\small $12$}
\cnode*(6,2){3pt}{N13}\rput(6,2.25){\small $13$}
\cnode*(5,1){3pt}{N14}\rput(4.75,1.25){\small $14$}
\cnode*(4,0){3pt}{N15}\rput(4,-0.35){\small $15$}
\cnode*(5,-1){3pt}{N16}\rput(4.75,-1.25){\small $16$}
\cnode*(-6,-4){3pt}{N17}\rput(-6,-4.35){\small $17$}
\cnode*(-5,-3){3pt}{N18}\rput(-5.25,-2.75){\small $18$}
\cnode*(-4,-2){3pt}{N19}\rput(-4.25,-1.75){\small $19$}
\cnode*(-4.5,-3.5){3pt}{N20}\rput(-4.15,-3.25){\small $20$}
\cnode*(-4.5,-4.5){3pt}{N21}\rput(-4.15,-4.75){\small $21$}
\cnode*(-5,-5){3pt}{N22}\rput(-5.25,-5.25){\small $22$}
\cnode*(-4,-6){3pt}{N23}\rput(-4.25,-6.25){\small $23$}
\cnode*(-3,-5){3pt}{N24}\rput(-2.75,-5.25){\small $24$}
\cnode*(-2,-4){3pt}{N25}\rput(-1.75,-3.75){\small $25$}
\cnode*(-3,-3){3pt}{N26}\rput(-2.75,-2.75){\small $26$}
\ncline{N1}{N2}
\ncline{N2}{N3}
\ncline{N4}{N5}
\ncline{N5}{N6}
\ncline[linewidth=2pt]{N2}{N5}\rput*(-0.45,-0.25){\small $e_3$}
\ncline{N7}{N8}
\ncline{N7}{N12}
\ncline{N8}{N9}
\ncline{N8}{N10}
\ncline{N9}{N16}
\ncline{N10}{N11}
\ncline{N10}{N14}
\ncline{N11}{N12}
\ncline{N11}{N16}
\ncline{N12}{N13}
\ncline{N13}{N14}
\ncline{N14}{N15}
\ncline{N15}{N16}
\ncangles[angleA=90,angleB=90,arm=3cm,linearc=.15]{N1}{N7}
\ncline{N3}{N13}
\ncline{N6}{N15}
\ncline{N4}{N9}
\ncline{N17}{N18}
\ncline{N17}{N22}
\ncline{N18}{N19}
\ncline{N18}{N20}
\ncline{N19}{N26}
\ncline{N20}{N21}
\ncline{N20}{N24}
\ncline{N21}{N22}
\ncline{N21}{N26}
\ncline{N22}{N23}
\ncline{N23}{N24}
\ncline{N24}{N25}
\ncline{N25}{N26}
\ncangle[angleA=135,angleB=135,arm=1cm,linearc=0.15]{N17}{N3}
\ncangle[angleA=0,angleB=270,arm=0.5cm,linearc=0.15]{N25}{N6}
\ncangle[angleA=0,angleB=270,arm=0.5cm,linearc=0.15]{N23}{N4}
\ncline{N1}{N19}
\end{pspicture}
\end{minipage}
}
\end{center}
\caption{Labels for $P_{26}$}
\label{P26Fig}
\end{figure}

\begin{prop}\label{P26}
Let $L$ be a copy of $P_{18}$ and $R$ be a copy of $P_{10}$.
Choose $e_2=xy$ to be one of the two bold--edges in $L$ and let $f,g \in H$ be a gadget--pair in $R$.
Then the dot product $L \cdot R$ gives rise to a new odd $2$--factored snark $P_{26}$.
Moreover, the only bold--edge of $P_{26}$ is $e_3$, the edge of $P_{26}$ identified with the edge $e_3$ of $P_{18}$ (cf. Figure~\ref{P26Fig}).
\end{prop}

\Prf
Applying the construction given by Theorem~\ref{OddDot} to the chosen bold--edge $e_2 \in L$ (cf. Lemma~\ref{P18}) and gadget--pair $f,g \in R$ (cf. Lemma~\ref{PetGadg}), we obtain that the graph $P_{26}$ is an odd 2--factored snark.

Fix the labelling on $V(P_{26})$ as in Figure~\ref{P26Fig}. To find all possible bold--edges in $P_{26}$,
we only need to verify Definition~\ref{boldedge}$(i)$, since we have just proved that $P_{26}$ is odd $2$--factored.

To this purpose, as in Proof of Lemma~\ref{P18}, we have implemented a program, with the software package MAGMA, and computed that
the graph $P_{26}$ has the dihedral group $D_4$ as automorphism group, its edge--orbits are eight and its vertex--orbits are seven.
For each representative $v$ of the vertex--orbits, we have determined all the $2$--factors of  $P_{26}-v$.
The only vertex, for which $P_{26}-v$ has only odd $2$--factors, is $v=2$ (lying in the vertex--orbit $\{2,5\}$).
Hence, the only bold--edge in $P_{26}$ is $e_3$, since there is an edge--orbit $E_0:=\{(2,5)\}$ of $P_{26}$, and its edge correspond to $e_3$ (c.f. Figure~\ref{P26Fig}).
\qed

\

\noindent {\sc Construction of $P_{34}$}

\

\begin{figure}[h]
\begin{center}
\resizebox{6cm}{!}{
\begin{minipage}[h]{12cm}
\begin{pspicture}(-5,-6)(8,6)
\psset{unit=1}
\cnode*(-2,0){3pt}{N1}\rput(-2.25,0){\small $1$}
\cnode*(0,-2){3pt}{N2}\rput(0,-1.75){\small $2$}
\cnode*(0,2){3pt}{N4}\rput(0,1.75){\small $4$}
\cnode*(2,0){3pt}{N3}\rput(1.75,0){\small $3$}
\cnode*(8,0){3pt}{N5}\rput(8,-0.35){\small $5$}
\cnode*(7,-1){3pt}{N6}\rput(7.15,-1.25){\small $6$}
\cnode*(6,-2){3pt}{N7}\rput(6,-2.25){\small $7$}
\cnode*(6.5,-0.5){3pt}{N8}\rput(6.15,-0.75){\small $8$}
\cnode*(6.5,0.5){3pt}{N9}\rput(6.15,0.75){\small $9$}
\cnode*(7,1){3pt}{N10}\rput(7.15,1.25){\small $10$}
\cnode*(6,2){3pt}{N11}\rput(6,2.25){\small $11$}
\cnode*(5,1){3pt}{N12}\rput(4.75,1.25){\small $12$}
\cnode*(4,0){3pt}{N13}\rput(4,-0.35){\small $13$}
\cnode*(5,-1){3pt}{N14}\rput(4.75,-1.25){\small $14$}
\cnode*(-6,-4){3pt}{N15}\rput(-6,-4.35){\small $15$}
\cnode*(-5,-3){3pt}{N16}\rput(-5.25,-2.75){\small $16$}
\cnode*(-4,-2){3pt}{N17}\rput(-4.25,-1.75){\small $17$}
\cnode*(-4.5,-3.5){3pt}{N18}\rput(-4.15,-3.25){\small $18$}
\cnode*(-4.5,-4.5){3pt}{N19}\rput(-4.15,-4.75){\small $19$}
\cnode*(-5,-5){3pt}{N20}\rput(-5.25,-5.25){\small $20$}
\cnode*(-4,-6){3pt}{N21}\rput(-4.25,-6.25){\small $21$}
\cnode*(-3,-5){3pt}{N22}\rput(-2.75,-5.25){\small $22$}
\cnode*(-2,-4){3pt}{N23}\rput(-1.75,-3.75){\small $23$}
\cnode*(-3,-3){3pt}{N24}\rput(-2.75,-2.75){\small $24$}
\cnode*(-4,2){3pt}{N25}\rput(-4.25,1.75){\small $25$}
\cnode*(-3,3){3pt}{N26}\rput(-2.75,2.75){\small $26$}
\cnode*(-2,4){3pt}{N27}\rput(-2,3.75){\small $27$}
\cnode*(-3.5,3.5){3pt}{N28}\rput(-3.85,3.25){\small $28$}
\cnode*(-3.5,4.5){3pt}{N29}\rput(-3.85,4.75){\small $29$}
\cnode*(-3,5){3pt}{N30}\rput(-2.75,5.25){\small $30$}
\cnode*(-4,6){3pt}{N31}\rput(-4.25,6.25){\small $31$}
\cnode*(-5,5){3pt}{N32}\rput(-5.25,5.25){\small $32$}
\cnode*(-6,4){3pt}{N33}\rput(-6.25,4.25){\small $33$}
\cnode*(-5,3){3pt}{N34}\rput(-5.25,2.75){\small $34$}
\ncline{N5}{N6}
\ncline{N5}{N10}
\ncline{N6}{N7}
\ncline{N6}{N8}
\ncline{N7}{N14}
\ncline{N8}{N9}
\ncline{N8}{N12}
\ncline{N9}{N10}
\ncline{N9}{N14}
\ncline{N10}{N11}
\ncline{N11}{N12}
\ncline{N12}{N13}
\ncline{N13}{N14}
\ncangles[angleA=90,angleB=90,arm=3cm,linearc=.15]{N1}{N5}
\ncline{N4}{N11}
\ncline{N3}{N13}
\ncline{N2}{N7}
\ncline{N15}{N16}
\ncline{N15}{N20}
\ncline{N16}{N17}
\ncline{N16}{N18}
\ncline{N17}{N24}
\ncline{N18}{N19}
\ncline{N18}{N22}
\ncline{N19}{N20}
\ncline{N19}{N24}
\ncline{N20}{N21}
\ncline{N21}{N22}
\ncline{N22}{N23}
\ncline{N23}{N24}
\ncangle[angleA=135,angleB=135,arm=1cm,linearc=0.15]{N15}{N4}
\ncangle[angleA=0,angleB=270,arm=0.5cm,linearc=0.15]{N23}{N3}
\ncangle[angleA=0,angleB=270,arm=0.5cm,linearc=0.15]{N21}{N2}
\ncline{N1}{N17}
\ncline{N25}{N26}
\ncline{N25}{N34}
\ncline{N26}{N27}
\ncline{N26}{N28}
\ncline{N27}{N30}
\ncline{N28}{N29}
\ncline{N28}{N32}
\ncline{N29}{N30}
\ncline{N29}{N34}
\ncline{N30}{N31}
\ncline{N31}{N32}
\ncline{N32}{N33}
\ncline{N33}{N34}
\ncangle[angleA=225,angleB=225,arm=1cm,linearc=0.15]{N33}{N2}
\ncangle[angleA=0,angleB=90,arm=0.5cm,linearc=0.15]{N27}{N3}
\ncangle[angleA=0,angleB=90,arm=0.5cm,linearc=0.15]{N31}{N4}
\ncline{N1}{N25}
\end{pspicture}
\end{minipage}
}
\end{center}

\caption{Labels for $P_{34}$}
\label{P34Fig}
\end{figure}
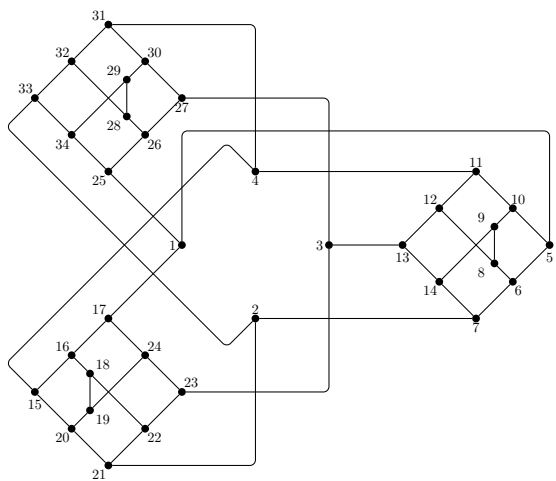

\begin{prop}\label{P34}
Let $L$ be a copy of $P_{26}$ and $R$ be a copy of $P_{10}$.
Choose $e_3=xy$ to be the only bold--edge in $L$ and let $f,g \in H$ be a gadget--pair in $R$.
Then the dot product $L \cdot R$ gives rise to a new odd $2$--factored snark $P_{34}$.
Moreover, $P_{34}$ has no bold--edges (cf. Figure~\ref{P34Fig}).
\end{prop}

\Prf
Applying the construction given by Theorem~\ref{OddDot} to the only bold--edge $e_3 \in L$ (cf. Proposition~\ref{P26})
and gadget--pair $f,g \in R$ (cf. Lemma~\ref{PetGadg}), we obtain that the graph $P_{34}$ is an odd 2--factored snark.

Fix the labelling on $V(P_{34})$ as in Figure~\ref{P34Fig}. To find all possible bold--edges in $P_{34}$, again,
we only need to verify that Definition~\ref{boldedge}$(i)$ does not hold, since we have just proved that $P_{34}$ is odd $2$--factored.

To this purpose, as in Lemma~\ref{P18} and Proposition~\ref{P26}, we have implemented a program, with the software package MAGMA, and computed that
the graph $P_{34}$ has the symmetric group $S_4$ as automorphism group, its edge--orbits and its vertex--orbits are both four.
For each representative $v$ of the vertex--orbits, we have determined all the $2$--factors of  $P_{34}-v$.
We have obtained that there is always a $2$--factor containing a cycle of even length. Thus, Definition~\ref{boldedge}$(i)$ does not hold. Hence $P_{34}$ has no bold--edges.
\qed

\

\begin{remark}\label{hagglund}
We have learned from J. H\"{a}gglund~\cite{Hag} that Brimnkmann, Goedbgebeur, Markstrom and himself had also found in~\cite{BGHM}, with an exhaustive computer search of all snarks of order $n\le 36$, numerical counterexamples to Conjecture~\ref{con1}, one of order $26$ and one of order $34$, but at the time we have informed him that we had already constructed these counterexamples via the bold--gadget dot product.
Indeed, we have checked that the snarks $P_{26}$ and $P_{34}$ are isomorphic to their graphs of order $26$ and $34$, respectively.
\end{remark}

\section{A partial characterization of odd $2$--factored snarks}\label{partcharact}


To approach the problem of characterizing all odd $2$--factored snarks, we consider
the possibility of constructing further odd $2$--factored snarks with the technique presented in Section~\ref{constrnewodd}, which
relies in finding other snarks with bold--edges and/or gadget--pairs,
Therefore, we study the existence of bold--edges and gadget--pairs in the known odd $2$--factored snarks.

We have already computed all the bold--edges in the Petersen graph $P_{10}$, the Blanu\v{s}a2 snark $P_{18}$, and the new snarks $P_{26}$ and $P_{34}$
(cf. Lemma~\ref{PetBold}, Lemma~\ref{P18}, Proposition~\ref{P26}, Proposition~\ref{P34}).

\begin{lemma}\label{FlowerNoBold}
Let  $J(t)$, for odd $t \ge 5$, be the Flower Snark. Then $J(t)$ has no bold--edges.
\end{lemma}

\Prf
Fix the labelling on the vertices of $J(t)$ as defined in Section~\ref{Pre}.
The flower snark has the dihedral group $D_{2t}$ as automorphism group~\cite{FR08},
its edge--orbits are four and its vertex--orbits are three.

To prove that there are no bold--edges, we only need to verify Definition~\ref{boldedge}$(i)$ does not hold, since we have already proved in~\cite{ALS}
that $J(t)$ is odd $2$--factored. To this purpose, we have to find a $2$--factor containing an even cycle in $J(t) - v$,
for each representative $v$ of the vertex--orbits.

Let $h_1,w_1$ and $u_1$ be representatives for the three vertex--orbits of $J(t)$.
Then, for each orbit we can construct the following $2$--factor in $J(t) - v$:

\noindent \resizebox{13cm}{!}{
\begin{minipage}[h]{15cm}
$$
\begin{array}{|c|l|c|}
  \hline
  \mbox{Graph} & \mbox{2--factor cycles} & \mbox{lengths} \\
  \hline
  J(t)-h_1 & (w_i,h_i,v_i,v_{i+1},h_{i+1},w_{i+1}) \mbox{ for } i=3,5,7, \ldots, t-2 & \frac{t-3}{2} \mbox{ cycles of length } 6 \\
           & (u_1,u_2, \ldots, u_t,v_1,v_2,h_2,w_2,w_1,w_t,h_t,v_t)              & \mbox{ a cycle of length } t+8 \\
  \hline
  J(t)-w_1 & (w_i,h_i,v_i,v_{i+1},h_{i+1},w_{i+1}) \mbox{ for } i=2,4, \ldots, t-1 & \frac{t-1}{2} \mbox{ cycles of length } 6 \\
            & (h_1,u_1,u_2, \ldots, u_t,v_1)                                     & \mbox{ a cycle of length } t+2 \\
  \hline
  J(t)-v_1 & (w_i,h_i,v_i,v_{i+1},h_{i+1},w_{i+1}) \mbox{ for } i=4,6, \ldots, t-1 & \frac{t-3}{2} \mbox{ cycles of length } 6 \\
           & (v_1,h_1,w_1,w_2,w_3,h_3,v_3,v_2,h_2,u_2,u_3, \ldots, u_t)          & \mbox{ a cycle of length } t+8 \\
  \hline
\end{array}
$$
\end{minipage}
}

\

Hence, we have obtained that, for all of these graphs, there is always a $2$--factor containing an even cycle. Thus, Definition~\ref{boldedge}$(i)$ does not hold. Hence, $J(t)$ has no bold--edges.
\qed

\

Regarding gadget--pairs, we have computed so far, only the gadget--pairs in the Petersen graph $P_{10}$ (cf. Lemma~\ref{PetGadg}).

\begin{lemma}\label{gadgetsearch}
Let the Flower snark $J(t)$, for odd $t \ge 5$, the Blanu\v{s}a2 snark $P_{18}$, $P_{26}$, and $P_{34}$ be defined as above.
Then
  \begin{enumerate}[(i)]
    \item $P_{18}$, $P_{26}$ and $P_{34}$ have no gadget--pairs;
    \item The Flower snark $J(t)$ has no gadget--pairs.
  \end{enumerate}
\end{lemma}

\Prf
For each of these graphs, we will verify that Definition~\ref{gadgetpair}$(i)$ or $(iv)$ does not hold.

$(i)$ Fix the labelling on $P_{18}$, $P_{26}$ and $P_{34}$ as in Figures~\ref{P18Fig}, \ref{P26Fig} and \ref{P34Fig}.
For these graphs, we have implemented a program, with the software package MAGMA, in which we compute
the edge--orbits under the action of the automorphism group; we consider all independent edges $g=cd$ from a chosen representative $f=ab$ of each edge--orbit and then we find all $2$--factors of $G-\{f,g\}$. If any such $2$--factors exist then condition \ref{gadgetpair}$(i)$ does not hold.
Otherwise, we choose one of the edges $\{ac,ad,bc,bd\}$ (cf. Definition~\ref{gadgetpair}$(iv)$), say $ac$, then
we compute all $2$--factors of $G-\{f,g\}+\{ac\}$ and, in each case, we find a $2$--factor
for which condition \ref{gadgetpair}$(iv)$ does not hold.

In the graphs $P_{18}$, $P_{26}$ and $P_{34}$ for each representative $f=ab$ of one of the edge--orbits,
there are several possible independent edges $g=cd$.

        $$\begin{array}{|c|c|c|c|}
          \hline
                                               & P_{18} & P_{26}   & P_{34} \\
          \hline
           \mbox{\small Number of edge--orbits} &   6    &     8    &    4   \\
          \hline
           \begin{minipage}[h]{6cm}
           \begin{center}
           \small Number of independent edges \\ for each edge--orbit representative
           \end{center}
           \end{minipage} & 22 & 34 & 46 \\
          \hline
          \end{array}$$

For most pairs there exists a $2$--factor of $G-\{f,g\}$, thus Condition \ref{gadgetpair}$(i)$ does not hold, whereas
the pairs $f,g$ for which $G-\{f,g\}$ has no $2$--factors are:
        $$\begin{array}{|c|c|c|c|c|c|c|c|c|}
          \hline
          \multicolumn{9}{|c|}{P_{18}} \\
          \hline
           f & (1,2) & (9,1)   & \multicolumn{2}{c|}{(2,6)} & \multicolumn{4}{c|}{(12,13)}    \\
          \hline
           g & (7,8) & (12,13) & (4,8) & (12,13)             & (4,8) & (11,3) & (15,5) & (17,7) \\
          \hline
          \end{array}$$
\begin{center}
\resizebox{13cm}{!}{
\begin{minipage}[h]{15cm}
          $$\begin{array}{|c|c|c|c|c|c|c|c|c|c|}
          \hline
          \multicolumn{10}{|c|}{P_{26}}\\
          \hline
           f & \multicolumn{2}{c|}{(2,5)} & (7,8)   & (7,12) & \multicolumn{5}{c|}{(10,11)}    \\
          \hline
           g & (10,11) & (20,21)          & (13,14) & (9,16) & (1,7) & (3,13) & (4,9) & (6,15) & (20,21) \\
          \hline
          \end{array}$$
\end{minipage}
} 
\end{center}

          $$\begin{array}{|c|c|c|c|c|c|c|c|}
          \hline
          \multicolumn{8}{|c|}{P_{34}}\\
          \hline
           f & (5,6)   & \multicolumn{6}{c|}{(8,9)} \\
          \hline
           g & (11,12) & (1,5) & (2,7) & (3,13) & (4,11) & (18,19) & (28,29) \\
          \hline
          \end{array}$$

            For each of these pairs of edges, $G-\{f,g\}+\{ac\}$ admits a $2$--factor $F$ in which the cycle $C_{ac}$
            has odd length, or $F$ has other even cycles besides $C_{ac}$, contradicting \ref{gadgetpair}$(iv)$.
            Hence, $P_{18}$, $P_{26}$ and $P_{34}$ have no gadget--pairs, since Definition~\ref{gadgetpair}$(i)$ or $(iv)$ does not hold.

$(ii)$ For the graphs $J(t)$, $t \ge 5$ odd, fix the labelling on the vertices of $J(t)$ as defined in Section~\ref{Pre}.

Recall that in a cubic graph $G$, a $2$--factor, $F$, determines a corresponding $1$--factor, namely $E(G)-F$.
In studying $2$--factors in $J(t)$ it is more convenient to consider the structure of $1$--factors.

           If $L$ is a $1$--factor of $J(t)$ each of the $t$ links of $J(t)$ contain precisely one edge from $L.$
           This follows from the argument in~\cite[Lemma 4.7]{AAFJLS}. Then, a $1$--factor $L$ may be completely specified
           by the ordered $t$--tuple $(a_1, \, a_2 , \, \ldots , \, a_t)$ where $a_i \in \{u_i, \, v_i, \, w_i\}$ for each
           $i=1,2, \ldots , \, t$ and indicates which edge in $L_i$ belongs to $L.$ Together these edges leave a unique spoke
           in each $IC_i$ to cover its hub. Note that $a_i \neq a_{i+1}$, $i=1,2, \ldots , \, t$ (i.e. they lie
           in different channels, for example if $a_i=u_i$, then $a_{i+1} \ne u_{i+1}$).
           To read off the corresponding $2$--factor $F$ simply start at a vertex in a base cycle at the first interchange.
           If the corresponding channel to the next interchange is not banned by $L$, proceed along the channel to the next
           interchange. If the channel is banned, proceed via a spoke to the hub (this spoke cannot be in $L$) and then along
           the remaining unbanned spoke and continue along the now unbanned channel ahead. Continue until reaching a vertex already
           encountered, so completing a cycle $C_1.$ At each interchange $C_1$ contains either $1$ or $3$ vertices.
           Furthermore as $C_1$ is constructed iteratively, the cycle $C_1$ is only completed when the first interchange is revisited.
           Since $C_1$ uses either $1$ or $3$  vertices from $IC_1$ it can revisit either once or twice. If $C_1$ revisits twice then
           $C_1$ is a hamiltonian cycle which is not the case. Hence it follows that $F$ consists of two cycles $C_1$ and $C_2.$

           Let $f,g$ be independent edges in $J(t)$. Since each of the $t$ links of $J(t)$ contain precisely one edge from any given $1$--factor $L$ of $J(t)$, each $2$--factor of $J(t)$ must contain exactly two edges of each link $L_i$. Therefore, if $f,g \in L_i$, for some $i \in \{1, \ldots, t\}$, then there is no $2$--factor of $J(t)$ avoiding both, i.e. Definition~\ref{gadgetpair}$(i)$ holds. Hence, to prove statement $(ii)$ in this case, we need to verify that Definition~\ref{gadgetpair}$(iv)$ does not hold. We need first to prove that for all other independent pairs that Definition~\ref{gadgetpair}$(i)$ does not hold, namely that $J(t)-\{f,g\}$ contains a $2$--factor.
           To this purpose we will define a $1$--factor $L$ of $J(t)$ containing both $f$ and $g$, thus giving rise to a $2$--factor $J(t)-L$ of $J(t)-\{f,g\}$.
           As noted above, the $1$--factors of $J(t)$ can be specified by an ordered $t$--tuple $(a_1, \ldots, a_t)$ with $a_i \in \{u_i,v_i,w_i\}$ for
           $i=1, \ldots, t$. We need to consider the following four cases:

           {\sc Case 1:} $f,g$ belong to different links, i.e. $f \in L_i$ and $g \in L_j$, with $i, j \in \{1, \ldots, t\}$, $i \ne j$.
           Suppose $f=b_ib_{i+1}$ and $g=c_jc_{j+1}$, for $i \ne j$. Choose any $t$--tuple $(a_1, \, a_2 , \, \ldots , \, a_t)$
           such that $a_i=b_i$ and $a_j=c_j$. Define $L$ to be the $1$--factor of $J(t)$ corresponding to $(a_1, \, a_2 , \, \ldots , \, a_t)$.
           Note that in the case $j=i+1$, $b_{i+1} \ne c_j$, since $f$ and $g$ are independent.

           {\sc Case 2:} $f$ belongs to a link and $g$ is a spoke of the same index, i.e. $f \in L_i$ and $g \in IC_i$ for some $i \in \{1, \ldots, t\}$.
           Suppose $f=b_ib_{i+1}$ and $g=h_ic_i$.
           Choose any $t$--tuple $(a_1, \, a_2 , \, \ldots , \, a_t)$
           such that $a_{i-1} \ne c_{i-1}$ and $a_i=b_i$.
           Define $L$ to be the $1$--factor of $J(t)$ corresponding to $(a_1, \, a_2 , \, \ldots , \, a_t)$.
           Since $f$ and $g$ are independent, $c_i \ne b_i$.

           {\sc Case 3:} $f$ belongs to a link and $g$ is a spoke of different index, i.e. $f \in L_i$ and $g \in IC_j$,
           with $i, j \in \{1, \ldots, t\}$, $i \ne j$.
           Suppose $f=b_ib_{i+1}$ and $g=h_jc_j$, for $i \ne j$.
           Choose any $t$--tuple $(a_1, \, a_2 , \, \ldots , \, a_t)$
           such that $a_i=b_i$, $a_{j-1} \ne c_{j-1}$, $a_j \ne c_j$.
           Define $L$ to be the $1$--factor of $J(t)$ corresponding to $(a_1, \, a_2 , \, \ldots , \, a_t)$.
           Moreover, choose $a_{j-1}$ and $a_j$ in different channels, which is always possible since there are
           three channels at each link and only one needs to be avoided.

           {\sc Case 4:} $f,g$ are both spokes in different interchanges, i.e. $f \in IC_i$ and $g \in IC_j$,
           for $i, j \in \{1, \ldots, t\}$, $i \ne j$.
           Note that $f,g$ cannot be two spokes in the same interchange since they are independent edges.
           Suppose $f=h_ib_i$ and $g=h_jc_j$, for $i \ne j$.
           Choose any $t$--tuple $(a_1, \, a_2 , \, \ldots , \, a_t)$
           such that $a_{i-1} \ne b_{i-1}$, $a_i \ne b_i$, $a_{j-1} \ne c_{j-1}$, $a_j \ne c_j$,
           which is always possible since there are three channels at each link.
           Define $L$ to be the $1$--factor of $J(t)$ corresponding to $(a_1, \, a_2 , \, \ldots , \, a_t)$.

\

In each case, the $2$--factor $F$ corresponding to $J(t)-L$ is well defined and it avoids both $f$ and $g$, thus Definition~\ref{gadgetpair}$(i)$
does not hold.

     \


This leaves us to prove that in the case $f,g \in L_i$ for some $i \in \{1, \ldots, t\}$, in which
Definition~\ref{gadgetpair}$(i)$ holds, Definition~\ref{gadgetpair}$(iv)$ does not hold.
To this purpose, we choose one of the edges $\{ac,ad,bc,bd\}$ (cf. Definition~\ref{gadgetpair}$(iv)$), say $ac$,
and find a $2$--factor for which Definition~\ref{gadgetpair}$(iv)$ does not hold.
Suppose, that $f=ab=a_ia_{i+1}$ and $g=cd=c_ic_{i+1}$, with $a_i,c_i \in \{u_i,v_i,w_i\}$,
and consider the graph $J(t)-\{f,g\}+\{ac\}=J(t)-\{f,g\}+\{a_ic_i\}$.

Recall that the flower snark has the dihedral group $D_{2t}$ as automorphism group (\cite{FR08}) with
vertex orbits $[h_1]:=\{h_i \, :\, i = 1, \ldots, t\}$, $[w_1]:=\{w_i \, :\, i = 1, \ldots, t\}$, and
$[u_1]:=\{u_i,v_i \, :\, i = 1, \ldots, t\}$. Then, w.l.o.g. we can consider the following two cases:

           {\sc Case a:} $a_i,c_i \in [u_1]$, say $a_i=u_1$ and $c_i=v_1$.

           In this case the graph $J(t)-\{f,g\}+\{u_1v_1\}$ admits a $2$--factor $F$ of type $[3, t, 6, \ldots, 6]$ with cycles
           $(u_1h_1v_1)$, $(w_1w_2 \ldots w_t)$, and $(u_ih_iv_iv_{i+1}h_{i+1}u_{i+1})$, for $i=2,4, \ldots, t-1$.

           {\sc Case b:}  $a_i \in [u_1]$ and $c_i \in [w_1]$, say $a_i=u_1$ and $c_i=w_1$ respectively.

           In this case the graph $J(t)-\{f,g\}+u_1w_1$ admits a $2$--factor $F$ of type $[3, t+6, 6, \ldots, 6]$ with cycles
           $(u_1h_1w_1)$, $(v_1v_2 \ldots v_th_tw_tw_{t-1}h_{t-1}u_{t-1}u_t)$, and
           $(u_ih_iw_iw_{i+1}h_{i+1}u_{i+1})$, for $i=2,4, \ldots, t-3$.

In both cases, the cycle of $F$ containing $a_ic_i$ is odd (of length $3$) and it has some even cycles as well,
implying that Definition~\ref{gadgetpair}$(iv)$ does not hold.

Therefore, we can conclude that $J(t)$ has no gadget--pairs.
\qed

\

The results obtained so far give rise to the following partial characterization:

\begin{theorem}\label{partialodd}
Let $G$ be an odd $2$--factored snark of cyclic edge--connectivity four that can be constructed from the Petersen graph and the Flower snarks using the bold--gadget dot product construction. Then $G \in \{P_{18}$, $P_{26}$, $P_{34}\}$.
\end{theorem}

\Prf
There is no possibility to construct other odd $2$--factored snarks from the Flower snarks $J(t)$, $t\ge 5$ odd, with the bold--gadget dot product construction by Lemma~\ref{FlowerNoBold} and Lemma~\ref{gadgetsearch}$(ii)$.

The Blanu\v{s}a2 snark $P_{18}$, $P_{26}$, and $P_{34}$ have been constructed iteratively via the bold--gadget dot product from the Petersen graph
using the existence of bold--edges in Petersen (Lemma~\ref{PetBold}), $P_{18}$ (Lemma~\ref{P18}), $P_{26}$ (Proposition~\ref{P26}) and gadget--pairs in Petersen (Lemma~\ref{PetGadg}). Since $P_{34}$ has no bold--edges by Proposition~\ref{P34} and $P_{18}$, $P_{26}$, $P_{34}$ have no gadget--pairs by Lemma~\ref{gadgetsearch}$(i)$, there is no possibility to apply the bold--gadget dot product any further to these graphs.
\qed

\

\begin{conj}\label{newconj}
Let $G$ be a cyclically $5$--edge connected odd $2$--factored snark. Then $G$ is either the Petersen graph or the Flower snark $J(t)$, for odd $t\ge5$.
\end{conj}

\begin{remark}\label{remconj}

\noindent $(i)$ A minimal counterexample to Conjecture~\ref{newconj} must be a cyclically $5$--edge connected snark
of order at least $36$ (cf. Remark~\ref{hagglund}).
Moreover, as highlighted in~\cite{BGHM}, order $34$ is a turning point for several properties of snarks.

\noindent $(ii)$ It is very likely that, if such counterexample exists, it will arise from the superposition applied to one of the known odd $2$--factored snarks.

\noindent $(iii)$ We have also checked that the snark of order $46$, of perfect matching index $\tau(G)=5$, constructed by H\"{a}gglund in~\cite{Hag2}, counterexample to a strengthening of Fulkerson's conjecture~\cite{FV,Maz}, is not odd $2$--factored. Moreover, the Flower snark is odd $2$--factored but it has $\tau(G)=4$ (cf.~\cite{FV}).
Hence, there is no relation between odd $2$--factored snarks and their perfect matching index being $5$.

\end{remark}

\end{document}